\documentclass[11pt]{article}
\usepackage{amsfonts}
\newcommand{\C}{\mathbf{C}}
\newcommand{\bC}{\mathbf{\bar C}}
\renewcommand{\P}{\mathbf{P}^n}
\newcommand{\f}{\tilde{f}}
\def\eqnspace{\noalign{\vskip 2\jot}}
\def\Y{\mathfrak{Y}}
\def\O{\mathfrak{O}}
\def\dz{\fbox{dz}}
\def\ds{\displaystyle\strut}
\def\fb{\mathbf{f}}
\def\b{\mathbf{b}}
\def\Z{\mathbf{Z}}
\def\R{\mathbf{R}}
\def\U{\mathbf{U}}
\def\f{\tilde{f}}
\def\a{\mathbf{a}}
\def\dist{\mathrm{dist}}

\def\Lip{\mathrm{Lip}}
\def\Iso{\mathrm{Iso}}
\newtheorem{theorem}{Theorem}[section]
\newtheorem{corollary}[theorem]{Corollary}
\newtheorem{example}[theorem]{Example}
\newtheorem{proposition}[theorem]{Proposition}
\newtheorem{lemma}[theorem]{Lemma}
\title{Normal holomorphic curves from parabolic regions to
projective spaces}
\author{Alexandre Eremenko\thanks{Supported by NSF grant DMS-9800084}}
\begin{document}
\date{Spring 1998}
\maketitle
\begin{abstract}
A holomorphic map $\C\to\P$ is called normal if it is uniformly continuous
from the Euclidean metric
to the Fubini--Study metric. The paper contains a survey
of known results about such maps, as well as some new theorems.
\end{abstract}

\section{Holomorphic curves in projective spaces}

This text was written in 1998 as an answer to a question
asked by Misha Gromov. As he soon 
answered this question himself \cite{Gromov2},
this preprint was not intended for publication.
However the recent activity in the subject
\cite{Tsukamoto,Winkel} suggests that the survey
part of this paper, which occupies most of it,
might be of some use to the researchers
in this area. The result which was new in 1998
is in the Appendix, and a stronger result
is available now,
\cite[Thm. 1.5]{Tsukamoto}.
\vspace{.1in}

We consider holomorphic maps $f:G\to\P$, where $G$ is a region in
the complex line $\C$,
and $\P$ is the complex projective space
of dimension $n$. Such maps are
called holomorphic curves.

We denote by $\Pi:\C^{n+1}\backslash\{ 0\}\to\P$ the standard projection map.
If $w=\Pi(\zeta)$ and
$\zeta=(\zeta_0,\ldots,\zeta_n)$ we write
$w=(\zeta_0:\ldots:\zeta_n)$ (with columns) and call $\zeta_j$
homogeneous coordinates of $w$.
The Fubini--Study metric is given in homogeneous coordinates by
\begin{equation}
\label{metric}
ds^2=\frac{\langle d\zeta,d\zeta\rangle\langle
\zeta,\zeta\rangle-
|\langle\zeta,d\zeta\rangle|^2}{\langle \zeta,\zeta\rangle^2},
\end{equation}
where $\langle.\,,.\rangle$ stands for the standard Hermitian product
in $\C^{n+1}$. The Fubini--Study distance between two points $\Pi(\zeta)$
and $\Pi(\eta)$ is equal to the ``angle'' between two complex one-dimensional
subspaces passing through 
$\zeta$ and $\eta$ 
that is 
\begin{equation}
\label{distance}
\dist(\Pi(\zeta),\Pi(\eta))=\arccos\frac{|\langle\zeta,\eta\rangle|}{\| \zeta\|\|
\eta\|},\quad\mbox{where}\quad \| \zeta\|:=\sqrt{\langle\zeta,\zeta\rangle} 
\end{equation}
is the Euclidean norm in $\C^{n+1}$. So, for example, the diameter of $\P$
is equal to $\pi/2$.
In what follows all metric notions in
$\P$ will refer to the Fubini--Study metric, and in
$\C^k$ to the Euclidean metric. We use $\dist$ and $B(a,r)$ for distances
and open balls in both cases, in addition the Euclidean distance will 
sometimes be written using the norm notation $\| a-b\|$.

Every holomorphic curve can be
factored as $f=\Pi\,\circ\f$ where
$\f:G\to\C^{n+1}\backslash\{0\}$ is a holomorphic
map, called a {\em reduced homogeneous representation} of $f$. Thus
$\f=(f_0,\ldots,f_n)$, where $f_j$ are holomorphic functions without common zeros.
For a given curve $f$ its reduced homogeneous representation is defined 
up to multiplication of all coordinates by the same holomorphic function
without zeros.
We will also use {\em meromorphic homogeneous representations}, where
the homogeneous coordinates are allowed to be meromorphic functions and
to have common zeros. Every $(n+1)$-tuple of meromorphic functions 
$(f_0,f_1,\ldots,f_n)$ in $G\subset\C$ defines a holomorphic curve
$f:G\to\P$, except when all functions $f_j$ are identically equal to $0$.
Indeed, we can always multiply all coordinates by a meromorphic function
in $G$ to cancel all common zeros and all poles.
{\em Unless otherwise is explicitly stated we use everywhere only reduced
homogeneous representations.}
If $n=1$ the curve $f=(f_0:f_1)$ is identified
with the meromorphic function $f=f_1/f_0: \C\to\bC={\bf P}^1$. Here $\bC$
is the Riemann sphere; the Fubini-Study metric for $n=1$ is the
spherical metric of constant curvature~$4$.

The length distortion of a holomorphic curve $f$
(from the Euclidean to the Fubini--Study metric)
is described by the {\em spherical derivative} $f^\#$, and the area
distortion by its square.
The following explicit expression can be derived from (\ref{metric}): 
\begin{equation}
\label{derivative}
({f}^\#)^2:=\frac{\sum_{i<j}|f_i^\prime f_j-f_if_j^\prime|^2}{\|\f\|^4}.
\end{equation}
Here $\f=(f_1,\ldots:f_n)$ is a reduced homogeneous representation 
($f^\#$ does not depend on its choice).
When $n=1$ we have
$$f^\#=\frac{|f'|}{1+|f|^2}.$$

An introduction to Nevanlinna--Cartan theory
is \cite{Lang}.
The Nevanlinna--Cartan characteristic of a holomorphic curve $f:\C\to\P$
is defined by
\begin{equation}
\label{T}
T(r,f)=\frac{1}{2\pi}\int_{-\pi}^\pi\log\|\f(re^{i\theta})\|\,d\theta-
\log\|\f(0)\|,
\end{equation}
where $\f$ is a reduced homogeneous representation. It is easy to see
that $T$ does not depend on the choice of representation.
The Laplacian of the subharmonic function $\log\|\f\|$ has density
$2(f^\#)^2$ with respect to the Lebesgue measure $\dz$ in the plane, so by
Jensen's Formula we have
\begin{equation}
\label{as}
T(r,f)=\int_0^rA(t,f)\frac{dt}{t},\quad
\mbox{where}\quad A(t,f)=\frac{1}{\pi}\int_{|z|\leq t}(f^\#)^2(z)\dz. 
\end{equation}
If $n=1$ this is called the Ahlfors--Shimizu
form of the Nevanlinna characteristic,
and $A(r,f)$ Ahlfors' (non-integrated) characteristic. It
is equal to the area of the disc $B(0,t):=\{ z:|z|\leq t\}$ with
respect to the pull-back of the spherical metric, divided by $\pi$.
As the total area of the Riemann
sphere is equal to $\pi$, the non-integrated characteristic $A(t,f)$
can be interpreted as the average covering degree of $f:B(0,t)\to\bC$.

The order of a curve $f:\C\to\P$ is defined by
$$\rho_f:=\limsup_{r\rightarrow\infty}\frac{\log T(r,f)}{\log r}.$$
If $f$ is a curve of finite order $\rho$, there exists a reduced homogeneous
representation whose coordinates have order at most $\rho$.

If $f:\C^*\to\P$ then the definitions of characteristics have
to be slightly modified. We put
$$A(r,f)=\frac{1}{\pi}\int_{\{ z:0\leq\log|z|/\log r\leq 1\}}
(f^\#)^2(z)\dz,\quad r>0$$
and  
$$T(r,f)=\frac{1}{2\pi}\int_{-\pi}^\pi\log\|\f(re^{i\theta})\|d\theta-
\frac{1}{2\pi}\int_{-\pi}^\pi\log\|\f(e^{i\theta})\|d\theta.$$
Then we have again the first relation in (\ref{as}).
There are two values of the order now: $\rho_f(0)$ and $\rho_f(\infty)$,
one for each singularity.
%where $\theta\in[0,2\pi)$ is chosen so that

\section{Normal curves in parabolic regions}

The set of all holomorphic curves $G\to\P$, equipped
with topology of uniform convergence on compacts in $G$ with respect
to the Fubini--Study metric, forms a complete metric space. A set of
holomorphic curves in a region $G\in\C$ is called
a {\em normal
family} if the closure of this set is compact. A necessary and sufficient condition for normality is that
the family is equicontinuous on every compact subset of $G$
(Ascoli--Arzela Theorem). An equivalent way to say this is that
spherical derivatives are uniformly bounded on on compacts in $G$.

Every region $G\in\bC$ has a complete Riemannian metric of constant
curvature, compatible with the conformal structure. For a given region 
such metric is defined up to a constant multiple. We choose the following
normalizations.
In $\bC$ we take the 
spherical metric defined above, it has curvature $4$.
In $\C$ we choose the standard Euclidean metric and in $\C^*$ the Riemannian
metric $|dz|/|z|$, both of them of zero curvature.
In a hyperbolic region there exists unique complete conformal
metric of 
curvature $-4$ which comes from the metric $|dz|/(1-|z|^2)$ in the unit disc
$\U$ via the
Uniformization Theorem. We call these metrics intrinsic for $G$. 
The group of isometries  
is denoted by ${\Iso}(G)$. 
There are four regions, namely $\bC,\,\C,\,\C^*$ and $\U$ where the group of
isometries acts transitively.

Let $G$ be a region whose group of isometries acts transitively.
A holomorphic curve $f:G\to\P$ is called {\em normal} if it satisfies the 
following equivalent conditions
\begin{enumerate}
\item  
The family $\{ f\circ\phi:\phi\in{\Iso}(G)\}$ is normal.
\item 
$f$ is uniformly continuous from the intrinsic metric of $G$ to the Fubini--Study metric.
\item
$\sup_{z\in G}f^\#(z)/\rho(z)<\infty$, where $\rho$ is the ratio of
the intrinsic metric to the Euclidean metric. 
\end{enumerate}
\hfill$\Box$
\medskip

We reserve the name {\em normal function} for the case $n=1$.
The set of normal curves will be denoted by $\Y_{G,n}$, or $\Y_n$ if
$G=\C$.
For $K>0$ we set $\Y_{G,n}(K)=\{ f\in\Y_{G,n}:\sup\,(f^\#/\rho)\leq K\}$.
For every $G,K$ and $n$ the set $\Y_{G,n}(K)$ is compact and the group
${\Iso}(G)$ acts on it by translations:
$f\mapsto \phi f=f\circ\phi^{-1},\;\phi\in{\Iso(G)}.$

The subject of this paper is normal holomorphic curves defined in
parabolic regions $\C$ and $\C^*$. The elements of $\Y_1$
are called sometimes Yosida functions. They were introduced by Julia
\cite{Julia} and studied by Yosida \cite{Y}.
The importance of the class $\Y_n$ is
partially explained by the following theorem, based on the idea
of Lohwater and Pommerenke \cite{Lohwater} (case $n=1$). The same idea
was used effectively by Brody   
\cite[Ch. III]{Lang}, for holomorphic curves to compact manifolds.
See also \cite{Zalcman,Minda}.
\begin{theorem}
\label{rescaling}
Let $M$ be a set of holomorphic curves in $\C$, containing non-constant
curves and having the following properties:
\begin{enumerate}
\item[$(i)$]
if $f\in M$ and $L(z)=az+b,\; a\neq0$, 
then $f\circ L\in M$;
\item[$(ii)$]
$M\cup\{\mbox{\rm constant curves}\}$ is closed. 
\end{enumerate}
Then $M$ contains non-constant normal curves.
\end{theorem}
This theorem is useful because in some cases it permits to reduce
Picard-type theorems to their special cases for curves in $\Y_n$.
The requirement (i) can be substantially relaxed \cite{Pang1,Pang}.
Examples of applications 
are in \cite{Bergweiler,Bonk,Chen,Wong}; the survey of related results
in dimension $1$ is 
\cite{Zalcman1}. 

{\em Proof of Theorem \ref{rescaling}}. Let $f\in M$ be a non-constant curve. Put
$$M_n:=\max_{|z|\leq n}(n-|z|)f^\#(z):=(n-|z_n|)f^\#(z_n),\quad |z_n|<n.$$
Evidently $M_n\to\infty$. So
\begin{equation}
\label{star1}
\rho_n:=\frac{1}{f^\#(z_n)}=\frac{n-|z_n|}{M_n}=o(n-|z_n|).
\end{equation}
We put $g_n(z)=f(z_n+\rho_nz)$. Then
\begin{equation}
\label{11}
g_n^\#(0)=\rho_nf^\#(z_n)=1,
\end{equation}
and for any fixed $r>0$ and $|z|\leq r$ we have, using (\ref{star1})
and the definition
of $M_n$:
$$\displaystyle\begin{array}{l}\ds g_n^\#(z)=\rho_nf^\#(z_n+\rho_nz) \\
\eqnspace
\displaystyle\leq\frac{n-|z_n|}{M_n}\max_{|z|\leq|z_n|+\rho_nr}f^\#(z) \\
\eqnspace
\displaystyle\leq\frac{n-|z_n|}{n-|z_n|-\rho_n r}.\frac{1}{M_n}.\max_{|z|\leq|z_n|+\rho_nr}
(n-|z|)f^\#(z) \\
\eqnspace
\displaystyle\leq (1+o(1))\frac{1}{M_n}\max_{|z|\leq n}(n-|z|)f^\#(z)=(1+o(1)).\end{array}$$
Thus $\{g_n\}$ is a normal family, and we can choose a subsequence such that
$g_n\to g$, where $g$ is non-constant holomorphic curve in view of (\ref{11}).
The assumption (i) implies that $g_n\in M$, and thus by (ii) we have $g\in M$.
We also have $g^\#\leq 1$, so $g$ is normal. \hfill$\Box$

\medskip
 
It follows from (\ref{as}) that $f\in\Y_n(K)$ satisfy
\begin{equation}
\label{ordertwo}
T(r,f)\leq K^2r^2/2,\quad r>0,
\end{equation}
so they are {\em of order at most $2$, normal type}.

The following characterization of $\Y_n$ belongs to Montel and Yosida
for $n=1$.
A set of hypersurfaces in $\P$ is called admissible if every $n+1$ hypersurfaces
of this set have empty intersection.
\begin{theorem}
\label{3n+1} Let $H_1,\ldots,H_{3n+1}$ be an admissible set of
hypersurfaces and $f$ a holomorphic curve. Denote by $E_j=f^{-1}(H_j)$ the
preimages of these hypersurfaces. Then $f\in\Y_n$ if and only
if the following condition is satisfied: there exists $\delta>0$ such that
every disc of diameter $\delta$ in $\C$ intersects at most $n$ of the
sets $E_j$.
\end{theorem}

{\em Remark}. For $n=1$ we need $4$ points (any set of points in $\bC$ 
is admissible). There are examples showing that three points may not
be enough. 

{\em Proof of Theorem \ref{3n+1}}. Let $f\in\Y_n(K)$, 
and $B\subset\C$ a disc of diameter $\delta$. Then   
${\rm diam}f(B)\leq K\delta$, so if $\delta$ is small enough, $f(B)$
cannot intersect
$n+1$ hypersurfaces. 
Otherwise there would be a sequence of balls $B(w_k,r_k)\subset\P$ with
$r_k\to 0$, each ball intersecting $n+1$ hypersurfaces. By passing to
a subsequence we may assume that these $n+1$ hypersurfaces are
the same for all balls, say $H_1,\ldots,H_{n+1}$.
We can also assume that $w_k\to w\in\P$. But then
$w\in H_1\cap\ldots\cap H_{n+1}$ and this contradicts our assumption that
the system of hypersurfaces is admissible. This proves ``only if'' part
of Theorem~2.

Now we assume that for some $\delta>0$ every disc of diameter $\delta$
intersects at most $n$ of $E_j$. Fix such a disc $B$ and notice
that at least $2n+1$ hypersurfaces are omitted in $B$. It remains
to use the following generalization of Landau's theorem:
{\em if a holomorphic curve in a disc omits $2n+1$ hypersurfaces from
an admissible system then its spherical derivative is bounded on every
compact in this disc by a constant, depending only on the hypersurfaces and
the compact} \cite{Eremenko,Eremenko2}. \hfill$\Box$

\medskip

\begin{theorem}
\label{ratios}
Let $f=(f_0:\ldots:f_n)$ be a holomorphic curve.
If all ratios $f_i/f_j$ belong to $\Y_1(K)$
then $f\in\Y_n(K\sqrt{n})$.
\end{theorem}
The converse is not true as the following example shows:
$f(z)=(\cos z:\cos(\alpha z):z)$, where $\alpha\in(0,1)$ is irrational.
To show that $f\in\Y_2$ consider a sequence $\lambda_k\to\infty$.
By choosing a subsequence we may assume that
one of the following cases holds:
\begin{enumerate}
\item[(i)]
$(\cos\lambda_k)/\lambda_k\to\infty,$
\item[(ii)]
$(\cos\lambda_k)/\lambda_k\to 0,\quad\mbox{or}$
\item[(iii)]
$(\cos\lambda_k)/\lambda_k\to a\in\C^*.$
\end{enumerate}
It is easy to see that the translations $t_{\lambda_k}f$ converge
to $(1:0:0),\,(0:0:1)$
or $(a\exp(\pm 2iz):0:1)$ in cases (i),(ii) and (iii) respectively.
So $f\in\Y_2$. On the other hand, meromorphic function
$\cos z/\cos(\alpha z)$ does not belong to $\Y_1$ because some of its
zeros are very close to poles, which cannot happen for a uniformly continuous
function.

{\em Proof of Theorem \ref{ratios}}. If the spherical derivatives of all
ratios are at most $K$, we have
$$\sum_{i<j}|f_i^\prime f_j-f_if_j^\prime|^2\leq
K^2\sum_{i<j}(|f_i|^2+|f_j|^2)^2\leq K^2n(|f_1|^2+\ldots+|f_n|^2)^2.$$
\hfill$\Box$

\medskip

The following results about free interpolation for Yosida curves were stated 
by M. Gromov in his lecture in Tel Aviv University in November 1997.

A set $E\in \C$ is called $K$-{\em sparse} if the distance between points 
of $E$ is at least $K$.
\begin{theorem}
\label{existence}
There exists $C(n)>0$, depending only on dimension,
such that every function $E\to\P$,
defined on an $K$-sparse set $E$, can be interpolated by a $f\in\Y_n(C(n)/K)$.
\end{theorem}
The proof is given in the Appendix. It is not clear whether a similar
result is true with a constant $C$
independent on dimension.

A set $E\in\C$ is called $K$-{\em dense} if every square with side length $K$ contains
at least one point of $E$.
\begin{theorem}
\label{uniqueness}
For every $K$-dense set $E$ and every
$f_0$ and $f_1$ in $\Y_n(c/K)$ with  $c<\sqrt{\pi/2}$ the equality
$f_0|E=f_1|E$ implies $f_0=f_1$.
\end{theorem}

{\em Proof of Theorem \ref{uniqueness}}. Let $n_E(r)={\rm card}\{ z\in E:|z|\leq r\}.$
Then $n_E(r)\geq \pi (r/K)^2+O(r)$ so, assuming wlog that $0\notin E$,
\begin{equation}
\label{N}
N_E(r):=\int_0^rn_E(t)\frac{dt}{t}\geq \frac{\pi r^2}{2K^2}+O(r).
\end{equation}
On the other hand the assumptions of Theorem \ref{uniqueness}
 and (\ref{ordertwo}) imply
$$T(r,f_j)\leq (cr)^2/(2K^2),\quad j=0,1.$$
Now take a fractional-linear function $L:\P\to\bC$, and put
$g_j:=L\circ f_j,\;j=0,1.$ 
Then we have $T(r,g_j)\leq T(r,f_j)+O(1)$ and $T(r,g_0-g_1)\leq
T(r,g_0)+T(r,g_1)+
\log 2$. Thus $T(r,g_0-g_1)\leq (cr/K)^2+O(1).$ Now (\ref{N})
implies that $N(r,0,g_0-g_1)\geq\pi r^2/(2K^2)+O(r)$, so we obtain from
the First Main Theorem of Nevanlinna that  $g_0=g_1$. As this conclusion
is valid for every fractional linear function $L$ we conclude
that $f_0=f_1$. \hfill$\Box$ 

\section{Normal functions in $\C^*$}
In this section we consider the class $\Y_{\C^*,n}$ consisting of
holomorphic curves $\C^*\to\P$, with $\sup\,|z|f^\#(z)<\infty$.
It is naturally isomorphic to the subclass of $2\pi i$-periodic
curves in $\Y_{\C,n}$.
The multiplicative group
of $\C^*$ acts on $\C^*$ by isometries $z\mapsto\lambda z,\;\lambda\in\C^*$.
So our curves $f:\C^*\to\P$ are characterized by the property
that the families of their translations
$\{ h_\lambda:\lambda\in\C^*\}$, where $h_\lambda f(z):=f(\lambda z)$,
are normal. The class $\Y_{\C^*,1}$ will be also called $\O_1$ after 
A. M. Ostrowski, who gave in \cite{Ostr} an explicit parametric
description of the class
(see Theorem \ref{ostrowski} below). In fact 
Ostrowski considered slightly smaller class of functions which have
no essential singularity at $0$. This subclass was introduced by Julia
in connection with the so-called Julia directions. Ostrowski functions
with no essential singularity at $0$ are exactly meromorphic functions
in $\C$ which have no Julia directions. This subclass was studied by
Julia, Montel and Ostrowski under the name of ``exceptional functions''.
(``Exceptional'', because they have no Julia directions).
Chapter VI of \cite{Montel} contains a detailed exposition of this work,
including the remarkably complete result of Ostrowski, which we
slightly generalize here in
\begin{theorem}\label{ostrowski} A meromorphic function $f:\C^*\to\bC$
belongs to $\O_1(K)$ 
if and only if it admits a representation
\begin{equation}\label{repr}
f(z)=az^m
\frac{\prod_{k\geq 0}\left(1-\frac{z}{a_k}\right)
      \prod_{k<0}\left(1-\frac{a_k}{z}\right)}{
      \prod_{k\geq 0}\left(1-\frac{z}{b_k}\right)
      \prod_{k<0}\left(1-\frac{b_k}{z}\right)}, 
\end{equation}
where $a\in\C,\; m\in\Z,\; a_k\in\C^*,\; b_k\in\C^*$, both sequences $(a_k)$
and $(b_k)$
may be finite or infinite in one or both directions, they tend to
$0$ as $k\to-\infty$, and they tend to $\infty$ as $k\to+\infty$,
and the following four conditions are satisfied:
\begin{enumerate}

\item[$(i)$]
The  number of zeros $a_k$ and poles $b_k$ of $f$ 
in every ring of the form $\{ z:r<|z|<2r\},\;r>0$ (counting multiplicity)
is bounded by a constant $C_1(K)$.

\item[$(ii)$]
the difference between the number
of zeros  and poles of $f$ in every ring
$\{ z:r_1<|z|<r_2\},\; r_1>0,\; r_2>0$ (counting multiplicity)
is bounded by a constant $C_2(K)$.

\item[$(iii)$]
For every $p$ and $q$ the ratios
$$\displaystyle|a_p|^m\displaystyle\frac{\prod_{ k: 0\leq\log|a_k|/\log|a_p|\leq1}
\displaystyle\frac{|a_p|}{|a_k|}}{\prod_{ k: 0\leq\log|b_k|/\log|a_p|\leq 1}\displaystyle\frac{|a_p|}{|b_k|}}
\quad\mbox{and}\quad
|b_q|^m\displaystyle
\displaystyle\frac{\prod_{ k:0\leq\log|b_k|/\log|b_q|\leq1}
\displaystyle\frac{|b_q|}{|b_k|}}{\prod_{ k:0\leq\log|a_k|/\log|b_q|\leq 1}\displaystyle\frac{|b_q|}{|a_k|}}
$$
are bounded from above by a constant $C_3(K)$. 

\item[$(iv)$]
For every pair $(k,j)$ the distance between $a_k$ and $b_j$ is bounded
away from zero by a positive constant $C_4(K)$.
\end{enumerate}
\end{theorem}

There is a simple geometric interpretation of conditions (i)-(iii), also
given in
\cite{Ostr}.
The Jensen's formula for meromorphic functions in a ring $A(r_1,r_2)$
has the form
$$\frac{1}{2\pi}\int_{-\pi}^\pi\log |f(r_2e^{i\theta})|d\theta-
\frac{1}{2\pi}\int_{-\pi}^\pi\log |f(r_1e^{i\theta})|d\theta$$
$$=\int_{r_1}^{r_2}\left( n(t,0)-n(t,\infty)\right)\, d(\log t)+
s(\log r_2-\log r_1),$$
where $s$ is an integer. 
Now for a function $f$ of the form (\ref{repr}) we put 
$$\phi(t):=\frac{1}{2\pi}\int_{-\pi}^\pi\log|f(e^{t+i\theta})|d\theta.$$
Then $\phi$ is a piecewise linear function on $\R$. The jumps of derivative
correspond to zeros and poles of $f$. Namely,
the jump of derivative $\phi^\prime_+(t)-\phi^\prime_-(t)$ is equal to the
number of zeros minus the number of poles on the circle $\{ z:|z|=\exp(t)\}$.
So each time the derivative jumps by an integer.
Condition (i) implies that all jumps have bounded magnitude and
the number of jumps on any interval of length $\log 2$ is bounded.
Condition (ii) means that the algebraic sum of jumps on any interval
is bounded. This is equivalent to the boundedness of $\phi^\prime$ on the whole
real line.
Finally, condition (iii) means the following: there is a horizontal
strip such that whenever the graph of $\phi$ is above the strip, $\phi$ is
concave, and whenever the graph is below the strip, $\phi$ is convex.
Let us call piecewise linear functions with such properties {\em admissible}.

Once an admissible piecewise linear function is given,
one can construct function $f\in\O_1$ by prescribing the arguments
of zeros and poles at each point of jump of $\phi^\prime$, such that
condition (iv) is satisfied.

%Thus functions $f\in\O_1$ can be encoded by the following  objects.
%Consider a sequence of triples $(x_k,A_k,B_k)$, where the range of $k$
%is either all integers, or negative integers, or positive integers, $x_k\in\R$,
%and $A_k,B_k$ are finite subsets of the unit circle of bounded cardinality
%and with bounded distance between them. Denote by $|A|$ the 
%cardinality of a set. Consider the function 
%$$\phi_0(x)=\sum_{k:x_k\leq x}|A_k|-|B_k|$$
%belongs
%to

%OLet us call a {\em decorated broken line}
%a graph of a piecewise linear function, satisfying the above-mentioned
%properties, together with two sets of points on the unit circle assigned
%to each jump, one set representing arguments of poles another arguments of
%zeros on the corresponding circle. Then there is a one-to one correspondence
%between such decorated broken lines and functions $f\in\O_1$ up to a
%constant factor. This correspondence commutes with the action of $C^*$. 

Thus Theorem \ref{ostrowski} gives a simple and effective parametric
description of the class $\O_1$ in terms of their zeros, poles
and constants $a$ and $m$.  

%The following result shows that this parametrization is also continuous in
%a very strong sense. Denote by $Z(f)$ and $P(f)$ the sequences of zeros
%and poles of a function $f\in\O_1$ and by $a(f)$ and $q(f)$ the corresponding
%constants in the representation (\ref{repr}). Given two sequences
%$Z_1$ and $Z_2$ 
%(functions of subsets of integers) we define the distance in the following
%way: $\dist(Z_1,Z_2):=\inf
%\begin{theorem}
%For every $\epsilon>0$ and every $K>0$ there exists $\delta>0$ with the
%following property:
%\newline
%If $f_0$ and $f_1$ belong to $\O_1(K)$, and 

%Let $Z_i$, $P_i$ and $a_i$ denote the zero-set
%and the set of poles and the constant $a$ in (\ref{repr}) for $f_i$
%for $i=1,2$. Assume in addition that the natural integer $q$ in (\ref{repr})
%is the same for $f_0$ and $f_1$. 

Using (i) one can improve Theorem \ref{3n+1} for the case of functions
$\C^*\to\bC$: {\em Let $f$ be a meromorphic function in $\C^*$ and
$E_j=f^{-1}(a_j),\;j=1,2,3.$ Then for $f\in\O_1$ 
it is necessary and sufficient that
for some $\delta>0$ every disc of diameter $\delta$ in $\C^*$ intersects
at most one of the sets $E_j$.} See \cite[p. 162]{Montel}

It is interesting that for curves $\C^*\to\P$ there is a universal
lower bound for $\sup\,|z|f^\#(z)$. This was discovered by
Lehto and Virtanen \cite{LV} who used a geometric method;
then Lehto in \cite{Lehto} published
a very simple analytic proof with precise constant (for $n=1$).
We follow the method of
Lehto.
\begin{theorem}\label{lehto} For every non-constant holomorphic curve
$f:\C^*\to\P$ we have
$$\sup_{z}|z|f^\#(z)\geq 1/2.$$
So $\O_1(K)$ consists of constants for $K<1/2$.
\end{theorem}
%Let us apply this theorem to the case when $n=1$ and both
%singularities are removable. Let $g:\bC\to\bC$ be a topologically
%holomoprhic map (open and discrete). Then there exists a homeomorphism
%$\phi$, fixing $0$ and $\infty$ which makes $g\circ\phi$ holomorphic.
%So on this case Theorem (\ref{LV}) gives a pure geometrical
%result: for every open discrete map $g:\bC\to\bC$ and any STOP
{\em Proof}. There are two cases to consider.
\newline Case 1.
The order of $f$ at one of the singularities, $0$ or $\infty$
is positive. 
Then the statement of the theorem follows immediately from 
(\ref{as}).
\newline
Case 2. The order at both singularities is zero.
Let $f=(f_0:\ldots:f_n)$ be a reduced representation such that
all $f_j$ are functions of zero order at both singularities.
For every $w\in\C^*$ consider the function
$$g(z,w):=
f_0(z)\overline{f_0(\overline{zw})}+\ldots+f_n(z)\overline{f_n(\overline{zw})}.$$
Then $g$ is a holomorphic function in $\C^*\times\C^*$.
For fixed $w$ the function $g_w:z\to g(z,w)$ has zero order at
both singularities, so it either has zeros in $\C^*$ 
or has the form
\begin{equation}
\label{except}
g_w(z)=h(w)z^q,
\end{equation}
where $q$ is an integer.
If for some $w$ on the unit circle $g_w$ has a zero $z^*\in\C^*$
then the points $(f_0(z^*):\ldots:f_n(z^*))$
and $(f_0(\overline{z^*w}):\ldots: f_n(\overline{z^*w}))$ are ``diametrically
opposite points'' in $\P$ that is the distance between them
is equal to $\pi/2$. In this case $f$ assumes two diametrically
opposite values on the circle $|z|=|z^*|$. The intrinsic length
of this circle is equal to $2\pi$, so it follows that there
is a point $z$ on this circle where $|z|f^\#(z)\geq 1/2$.

It remains to consider the possibility that (\ref{except}) holds for 
all $w$ on the unit circle and all $z\in\C^*$. In this case
(\ref{except}) actually holds for all $z$ and $w$ in $\C^*$.
We examine this possibility by substituting for ${\tilde f}$
a Laurent series with undetermined coefficients:
$$f_j(z)=\sum_{k=-\infty}^{\infty}c_{j,k}z^k.$$
The resulting system of equations shows that the functional
equation (\ref{except}) has no solutions for which the curve
$f$ is not constant.\hfill$\Box$

\medskip

\noindent
{\bf Remarks}.

1. As curves $f\in \O_1$ have zero order, the second part
of the previous proof applies to them. It shows that for such
curve $f$ there always exist two points in the same circle
$\{ z:|z|=r\}$ whose images are diametrically opposite.
It is not clear whether such improvement of Theorem \ref{lehto}
is true for curves or order $\geq 1$.

2. The function $f(z)=z:\C^*\to\bC$
shows that the estimate $1/2$ in Theorem \ref{lehto} is precise.
Probably the only functions for which $\sup|z|f^\#(z)=1/2$ are
$f(z)= kz$ with 
$k\in\C^*$. The following example from \cite{Lehto} shows that for
every $\epsilon>0$ there is a ``periodic'' function $f:\C^*\to\bC$,
that is $f(tz)=f(z)$ for some $t>1$ and such that
$\sup|z|f^\#(z)\leq 1/2+\epsilon$. Put
$$f(z)=\prod_{k\in\Z}\frac{z+t^k}{z-t^k}.$$
Then direct computation shows that $\max_{\C*}|z|f^\#(z)\to 1/2$ as
$t\to\infty$.
This implies that there is an open discrete map from a flat torus whose
shortest closed geodesic has length $2\pi$ to the sphere whose
great circles have length $\pi$, such that the length distortion is
arbitrarily close to $1/2$. 
That a continuous map of non-zero degree with such properties exists
follows from Proposition 2.12 in\cite{Gromov}.

3. The proof in Case 2 shows that actually some circle $|z|=r$ has
image, not shorter than a great circle.
It is interesting to consider the special case when $f$ has no singularities,
that is extends to the whole sphere. For this case Theorem \ref{lehto}
follows from 
\begin{proposition}\label{gabrielov}
Let $f:\bC\to\P$ be a 
continuous map of non-zero degree. Then some circle $|z|=r$ has image of
length at least $\pi$.
\end{proposition}
This can be proved in the same way as Proposition 2.12 in \cite{Gromov}.

Our proof of the existence of free interpolation in the Appendix
shows that for periodic interpolation data one can find periodic
interpolating function. Thus Theorems \ref{existence} and
\ref{uniqueness} have counterparts for the class $\O_n$.

\section{Binormal curves}

A holomorphic curve $f\in\Y_{G,n}$ is called {\em binormal}
if the family of its translations by the elements of $\Iso(G)$ is
normal and has {\em no constant limit curves}.
Yosida in \cite{Y} called them functions of first category.
Notice that the closure of the set of translations of a binormal curve
consists of only binormal curves.
The following characterization was given by Yosida in \cite{Y} for $n=1$.
\begin{theorem}
\label{characterization}
A curve $f\in\Y_n$ is binormal if and only if for every
$\delta>0$ there exists a constant $c$ such that
$$\int_{|z-\zeta|<\delta}(f^\#)^2(z)\fbox{dz}\geq c\quad\mbox{for every}\quad
\zeta\in\C.$$
\end{theorem}
\hfill$\Box$
\begin{corollary}
Binormal curves $f$ satisfy
$$T(r,f)\asymp r^2.$$
\end{corollary}

\hfill$\Box$
The idea of the following result is contained in \cite{Y}, but the result
is stated there in a weaker form, and Yosida's
proof of it contains mistakes (his Lemma 1 is incorrect). We use the standard
notations of Nevanlinna theory, in particular,
for the averaged counting function $N_1$
of critical
points
\begin{theorem}
\label{N_1}
For binormal functions $f\in\Y_1$ we have
\begin{equation}
\label{ramextreme}
N_1(r,f)=2T(r,f)+O(1)
\end{equation}
and
\begin{equation}
\label{nodeficiencies}
N(r,a,f)=T(r,f)+O(1).
\end{equation}
\end{theorem} 
Both statements follow from
\begin{proposition}\label{bonk}
Let $f$ be a binormal function.
If the arc  $\alpha$ is the intersection of a circle
of radius greater than $1$
with a disc of radius $1$, then
\begin{equation}
\label{durak}
\int_\alpha |\log f^\#(z)|\,|dz|\leq c
\end{equation}
and for every $a\in\bC$
\begin{equation}
\label{eshcheodin}
\int_\alpha|\log \left(\dist(f(z),a)\right)^{-1}|\,|dz|\leq c
\end{equation}
where $c$ is a constant depending only on $f$.
\end{proposition}
See \cite[Lemma 2]{Y} or \cite{Bonk}.
Now (\ref{nodeficiencies}) is an immediate consequence
from (\ref{eshcheodin}) and the First Main Theorem of Nevanlinna,
and (\ref{ramextreme}) can be derived
from (\ref{durak}) as in \cite{Bonk}. It follows that the error term in
the Second Main Theorem of Nevanlinna is bounded for binormal functions $f$.
It is also proved in \cite{Bonk} that one can
``differentiate'' the asymptotic relations (\ref{ramextreme}) and
(\ref{nodeficiencies}). More precisely,
$$n_1(r,f)=(2+o(1))A(r,f),\quad r\not\in E,$$
and
$$n(r,a,f)=(1+o(1))(A(r,f)),\quad r\notin E,$$
where $A$ is the (non-integrated) Ahlfors' characteristic, $n_1$
is the number of critical points in $B(0,r)$, and the exceptional
set $E$ has zero density.

We have the following corollaries,
all of them due to Yosida.
\begin{corollary}\label{nodefs}
If $f\in\Y_1$ is binormal then  
it has no deficiencies, and even no Valiron deficiencies.
This implies that $f$ assumes every value in $\bC$.\hfill$\Box$
\end{corollary}
\begin{corollary}\label{preimage-dense}
If $f\in\Y_1$ is binormal then preimage of every point is
$K$-dense in $\C$ with some $K>0$.
\end{corollary}
Indeed if, say poles, are not dense there is a sequence of pole-free
squares with
centers $\lambda_k$ and side length tending to infinity.
As the family $t_{\lambda_k}f$ is normal, we can choose a subsequence
converging to $g$, where $g$ is also binormal.
On the other hand $g$ has no poles at all, which contradicts
Corollary \ref{nodefs}.\hfill$\Box$ 

\medskip

A property which follows immediately from the definition should be also
mentioned: {\em binormal functions have no asymptotic values}.

Thus binormal  functions display asymptotic
behavior, similar to elliptic functions. In certain sense such behavior
(described in preceding corollaries) is ``typical'' for meromorphic functions.
In the next section we consider the ``opposite extreme'' to (\ref{ramextreme}).

Using Theorem \ref{ostrowski} we can give explicit description
of binormal functions in $\O_1$, that are meromorphic
functions $f$ in $\C^*$ such that the families
$\{ h_\lambda f:\lambda\in\C^*\}$ are  normal without constant limit functions.
\begin{theorem}\label{ostrowski-binormal} A meromorphic function $f\in\O_1$ is  binormal
if and only if in 
representation $(\ref{repr})$ the following additional property
is satisfied:
\newline
$(v)$ There exists $C_5(f)>0$ such that every annulus
$\{ z:r<|z|<c_5r\},\; r>0$ contains at least one zero and at least 
one pole of $f$.
\end{theorem}
The additional property (v) implies that all ratios in (iii) are also bounded
away from zero or, which is equivalent,
the piecewise-linear
function $\phi$ introduced
after 
Theorem \ref{ostrowski} 
is bounded.

{\em Proof of Theorem \ref{ostrowski-binormal}}. The necessity of condition
(v) follows from Corollary \ref{preimage-dense}, applied to the binormal
function
$f\circ\exp\in\Y_1$. Sufficiency is evident,
because if condition (v) is satisfied, all limit functions have zeros and poles.

\hfill$\Box$

\section{Locally univalent and entire normal functions in~$\C$}
\begin{theorem}
The only locally univalent functions in $\Y_1$ are exponential functions
and fractional-linear functions. By exponential we mean
$L\circ\exp(az)$, where $L$ is a fractional linear transformation and
$a\in\C^*$.
\end{theorem}
This is proved in \cite{Bonk}. Using \cite{smalram} one can probably relax 
the assumption of local univalence in this theorem, by replacing it
with the weaker assumption $N_1(r,f)=o(T(r,f))$,
and still preserve the conclusion.

It follows from a theorem of Clunie and Hayman \cite{CH} that
for entire functions in $\Y_1$ the growth estimate (\ref{ordertwo})
can be substantially improved. We give a refined version of this
theorem based on the work of Pommerenke \cite{Pomm} and Minda \cite{Minda}.
\begin{theorem} \label{CHPM}
For an entire function $f$ the following conditions are equivalent:
\begin{enumerate}
\item[$(i)$]
$\sup\,f^\#\leq 1.$
\item[$(ii)$]
$|f(z)|\leq 1$ implies $|f'(z)|\leq 2$.
\item[$(iii)$]
$\nabla w\leq 2$, where $w:=\log^+|f|$.
\end{enumerate}
\end{theorem}
{\em Remarks}.
It follows from (iii) that $w(z)\leq 2|z|$ that is
$$|f(z)|\leq\exp\max\{2|z|,1\},$$
so $f$ has at most exponential type.

Let $E$ be a subset of $\C$. If $|E|\geq 5$ then the condition 
that $f^\#$ is bounded on $f^{-1}(E)$ is equivalent to $f\in\Y_1$ for
meromorphic function $f$. In the case of entire functions $|E|\geq 3$
is enough (see \cite{Hinkk} for these results). Thus the condition (ii)
in Theorem \ref{CHPM} can be replaced by much weaker condition.

{\em Proof}. The implication (i) $\Rightarrow$ (ii) and
(iii)$\Rightarrow$ (i) are evident.

To prove (ii)$\Rightarrow$ (iii) we set
$D=\{ z:|f(z)|>1\}$, and
for every $R>0$ consider the following function in $D$:
$$u_R=\frac{|f'|}{|f|(\log|f|+R)}.$$
Evidently
\begin{equation}
\label{22}
u_R(z)\leq 2/R,\quad z\in\partial D,
\end{equation}
and $u_R$ satisfies
\begin{equation}
\label{31}
\Delta\log u\geq u^2
\end{equation}
in the sense of distributions (by direct verification).
\begin{proposition}\label{prop}
If $u$ is a positive continuous function
in a plane domain $D$
with the properties $(\ref{22})$ and $(\ref{31})$ then $u\leq 2/R$ in $D$.
\end{proposition}

{\em Proof of Proposition \ref{prop}}.
Assume that our Proposition is not true, so $u(z_0)>2/R$ for some
$z_0\in D$. Consider the function 
$$v=\frac{2R}{R^2-|z-z_0|},\quad z\in B(z_0,R)$$
Evidently 
\begin{equation}
\label{44}
v(z)\geq 2/R,\quad z\in B(z_0,R)
\end{equation}
and (by direct verification)
\begin{equation}
\label{55}
\Delta\log v=v^2.
\end{equation}
Consider the set
$$K=\{ z\in D\cap B(z_0,R):u(z)>v(z)\}.$$
We have by assumption $u(z_0)>2/R=v(z_0)$ so $z_0\in K$.
Let $D_0$ be the component of $z_0$ in $K$. Then we have
$$u(z)=v(z),\; z\in\partial D_0$$ because the inequality $u(z)\leq v(z)$
holds for
$z\in\partial D$ (because of (\ref{22}) and (\ref{44})), and for 
$z\in\partial B(z_0,R)$ (because for such $z$ we have $v(z)=+\infty$).

On the other hand by (\ref{31}) and (\ref{55}) we have
$$\Delta(\log u-\log v)\geq u^2-v^2>0\quad\mbox{in}\; D_0$$
so $\log u-\log v$ is a positive subharmonic function in $D_0$, zero
on the boundary. This contradicts Maximum Principle, so the Proposition
is proved. \hfill $\Box$

\medskip

Applying Proposition \ref{prop} to our function $u_R$ we obtain
\begin{equation}
\label{star}
|f'|\leq\frac{2}{R}\max\{|f|,1\}(\log^+|f|+R)
\end{equation}
(in $D$ this follows from Proposition \ref{prop}, in the rest of the plane from
(ii). Because this is true for arbitrary $R>0$ we obtain
(by letting $R\rightarrow\infty$ with fixed $z$):
\begin{equation}
\label{7}
|f'|\leq 2\max\{|f|,1\}
\end{equation}
Now we put $w=\log^+|f|$, so that $\nabla w(z)=|f'/f|(z)$ when $|f(z)|>1$
and rewrite the last equation as
$\nabla w\leq 2,$ which is (iii).\hfill$\Box$.

\medskip

Notice that the auxiliary function $u_R$ used in the proof is
the density of the pull-back of the Poincare metric of the
punctured disc $\{ w:|w|>1/R\}$, and Proposition (\ref{prop}) is
a version of the Ahlfors--Schwarz lemma. 

\begin{example}\label{fryntov} There are entire functions in $\Y_1$
of every order
$\rho\in [0,1]$.
\end{example}
This example is due to A. Fryntov (private communication).
If $\rho=1$ we take $f(z)=\exp(z)$. Now assume that $\rho\in (0,1)$
put 
$$f(z)=\prod_{k=1}^\infty\left(1-\frac{z}{2^k}\right)^{[2^{k\rho}]},$$
where $[x]$ stands for the greatest integer $\leq x$. Subharmonic
unction $u:=\log|f|$
satisfies the approximate functional equation
$$u(2z)=2^\rho u(z)+u_0(z),$$
where $u_0$ is negligible. This permits to verify (iii) in Theorem \ref{CHPM}.
To construct an example of order $0$ we replace $\rho$ in the previous
formula by $\rho(k):=1/\log^+k$.\hfill$\Box$ 

\medskip

%An entire function $f$ is said to be of {\em bounded value distribution}
%if there exist $R>0$ and $C>0$ with the following
%property: for every $a\in\C$ and for every disc $D$ of
%radius $R$ the number of solutions of the equation $f(z)=a$ in $D$ is at
%most $C$ (counting multiplicity). Hayman \cite{H} gave the following
%characterization of this class:
%\begin{theorem}\label{bvd}
%An entire function $f$ has bounded value distribution if and only if
%there is a natural integer $N$ and $C>0$ with the property
%$$|f^{(k)}(z)|\leq C\max_{1\leq j\leq N}|f^{(j)}|(z)\quad\mbox{for every}
%\quad z\in\C.$$
%\begin{theorem} $BVD\subset Y$ (proper inclusion). 
%\end{theorem}
%This result does not extend to meromorphic functions (if meromorphic
%functions of bounded value distribution are defined in the same way as 
%entire ones).
%An example is $\cos z/\cos\alpha z$ with irrational $\alpha\in\R$.
%
%{\em Proof of Theorem \ref{bvd}}. Example \ref{fryntov} shows that
%there is an entire function in $\Y\backslash \BVD$.
%
%To prove the inclusion we consider
%the family $\{ f_\lambda:=\t_\lambda f:\lambda
%\in\C\}$, where $f\in\BVD$. This is a quasi-normal family in Montel's
%terminology \cite[p. 66]{Montel}, which means that from every sequence
%of $\lambda$'s(\lambda_k)$ we can extract a subsequence $(\lambda_k)$,
%such that $f_{\lambda_k}$ is convergent uniformly on compacts in $\C\backslash 
%E$, where $E$ is a discrete set. It is enough to show that in fact the
%sequence converges uniformly in a neighborhood of every point of $E$ as well.

%\section{Related classes of functions}

\section{Appendix. Proof of existence of free interpolation}

We start with establishing some notations. By an automorphism
of $\P$ we mean $\Pi\circ U\circ \Pi^{-1}$,
where $U$ is a unitary transformation
of $\C^{n+1}$. So automorphisms are biholomorphic isometries. 

Let $B(r)=\{ \zeta\in\C^n:\|\zeta\|<r\}$ be the open ball of radius $r$ centered
at the origin and ${\bar B}(r)$ its closure. Consider one of the standard
local coordinates in $\P$, namely 
\begin{equation}
\label{coord}
\psi:B(2)\to\P,\quad \psi(\zeta_1,\ldots\zeta_n)=(1:\zeta_1:\ldots:\zeta_n).
\end{equation}
The length distortion by $\psi$ is estimated by
\begin{equation}\label{ldn}
\frac{1}{1+\|\zeta\|^2}\leq\frac{ds}{\| d\zeta\|}\leq
\frac{1}{\sqrt{1+\|\zeta\|^2}}
\end{equation}
which implies 
\begin{equation}
\label{lipshitz}
\frac{1}{5}\leq\frac{ds}{\| d\zeta\|}\leq 1\quad\mbox{for}\quad
\zeta\in B(2).
\end{equation}
So $\psi^{-1}\in\Lip(5)$. 
Let $p_j$ be the automorphism of $\P$ which interchanges the homogeneous
coordinate number $0$ with the homogeneous coordinate number $j$, where
$j\in\{ 1,\ldots,n\}$. Let $V=B(2)\times\{0,\ldots,n\}$, let $p:V\to B(2)$
be the projection map, and $W=p^{-1}({\bar B}(1))\subset V$. We consider $V$
as a (disconnected) Riemannian manifold equipped with the pull-back of
the Euclidean metric via $p$. So the distance between different sheets
$B(2)\times\{ i\}$ and $B(2)\times\{ j\},\;j\neq i$, is infinite.
The following property is evident:
\begin{equation}
\label{complete}
\{ a\in V:\dist(a,W)\leq\epsilon\}\quad
\mbox{is a complete metric space if}\quad\epsilon<1.
\end{equation}
Evidently $p$ has a continuous right inverse
in every ball of radius $1$ centered at 
a point of ${\bar B}(1)$, and this inverse is an isometry onto the image.
A point $a\in V$ will be denoted by $a=(p(a);j)=(\zeta;j)
=(\zeta_1,\ldots,\zeta_n;j)$,
where $\zeta\in B(2)$ and $j\in\{ 0,\ldots,n\}$.
We use the notation $\| a^\prime-a^{\prime\prime}\|$ for the distance
between two points in $V$.

We define
$$\Psi:V\to\P,\quad \Psi(a)=\Psi(p(a);j)=p_j\circ\psi\circ p(a).$$
Then $\Psi$ is a surjective local diffeomorphism; in fact
\begin{equation} 
\label{surjective}
\Psi|W:W\to\P \quad\mbox{is surjective}.
\end{equation}
Indeed, for every point in $\P$ one can choose homogeneous coordinates
such that one coordinate is equal to $1$ and the rest have absolute
value at most~$1$.
Moreover,
\begin{equation}
\label{branch}
\begin{array}{l}
\mbox{for every point}\;\; a\in W\;\;\mbox{there exists right inverse}\\
\Psi^{-1}\;\;\mbox{defined in}\;\; B(p(a),\delta)\subset\P\;\;\mbox{with}\;\;
\delta:=1/11,\\
\mbox{such that}\;\;\Psi^{-1}(p(a))=a\;\;
\mbox{and}\;\;\Psi^{-1}\in\Lip(5).
\end{array}
\end{equation}
Now we can solve a two-point interpolation problem with uniform estimates.
\begin{lemma}\label{twopoint}
For every point $q\in \P$ there exists a map $g:\C\times V\to\P$,
with the following properties:
\begin{enumerate}

\item[$(a)$]
$z\mapsto g(z,a),\;\C\to\P$ is holomorphic and its spherical derivative
is uniformly bounded with respect to $a$.

\item[$(b)$]
The map $a\mapsto g(0,a):V\to\P$ has surjective restriction on $W$ and
satisfies the condition $(\ref{branch})$.

\item[$(c)$]
$\dist (g(z,a),q)\leq (25/4)|z|^{-3}<(1/4)$ for $|z|\geq 3$.

\item[$(d)$]
$\dist (g(z,a^\prime),g(z,a^{\prime\prime}))\leq
(5/4)\| a^\prime-a^{\prime\prime}\|\,|z|^{-4}$ for $|z|\geq 3$.
\end{enumerate}
\end{lemma}
{\em Proof}. We take in this proof
for convenience $q=(1:1:\ldots:1)$. Then the general
case can be obtained by composing $g$ with an appropriate automorphism of $\P$.
We will construct $g$ with properties (a), (c), (d) and $g(0,a)=\Psi(a),\; 
a\in V$.
Then (b) will be satisfied in view of (\ref{surjective}) and (\ref{branch}).
Let us assume for simplicity that $a=(\zeta_1,\ldots,\zeta_n;0)$. Construction
for $a=(\zeta;j)$ is then obtained by composition with $p_j,\; 1\leq j\leq n$.
(Notice that $p_j(q)=q$).
We put
$$g(z,a):=\left( g_0(z,a):\ldots:g_n(z,a)\right):=
(z^4+1:z^4+4z+\zeta_1:z^4+\zeta_2:\ldots:z^4+\zeta_n),$$
that is $g_j(z,a)=z^4+\zeta_j,\; 2\leq j\leq n$.
First we notice that for $|z|\geq 3$ the following estimates hold (we use
$|\zeta_j|<2$):
\begin{equation}
\label{approx}
|1-z^{-4}g_j(z,a)|\leq 5|z|^{-3}<1/5 
\end{equation}
\begin{equation}
\label{lip}
|g_j(z,a^\prime)-g_j(z,a^{\prime\prime})|\leq \|
a^\prime-a^{\prime\prime}\|
\end{equation}
Now (c) is easy to verify, using (\ref{metric}) and (\ref{approx}),
and (d) follows from (\ref{metric}), (\ref{approx}) and (\ref{lip}). 
%$$\arccos\frac{\left|\sum_{j=0}^n(1+\alpha_j)\right|}{
%\sqrt{n+1}\sqrt{\sum_{j=0}^n|1+\alpha_j|^2}}\leq\sqrt{5}\alpha,$$
%where $\alpha:=\max_j|\alpha_j|<1/5$.

It remains to prove (a).
{}From (c) follows that $g^\#$ is uniformly bounded with respect to
$a$ for $|z|>4$ (Cauchy estimate for derivatives).
 Now let $|z|\leq 4$. It is enough to show that
the function $(z,a)\to g^\#(z,a)$ is continuous for $a=(\zeta;0),
\;(z,\zeta)\in \bar{B}(4)\times {\bar B}(2)$.
In view of (\ref{derivative}) it is enough to check that
the denominator in (\ref{derivative}) is never equal to zero for
$(z,\zeta)\in \bar{B}(4)\times \bar{B}(2)$, that is $g_j$ never have
common zero.
In fact we will show that the first two coordinates $g_0$ and $g_1$ never
have common zero. The only zeros of $g_0$ are $\pm1,\pm i$. So our assertion
follows from the fact that $g_1$ never has zeros in the ring
$\{ z: 3/4<|z|<5/4\}$ which follows from Rouche theorem.\hfill$\Box$

\medskip

Functions $g$ interpolate at two points, $0$ and $\infty$.
Now we are going to combine our two-point solutions, essentially by adding them,
so we need a surrogate of addition in $\P$.
Let 
\begin{equation}
\label{q}
\Delta:=B(\delta)=B(1/11)\subset\C^n\quad\mbox{and}\quad q:=(1:0:\ldots:0)\in\P. 
\end{equation}
These notations will be fixed till the end of the proof of
Theorem \ref{existence}.
We define a map $P:\P\times\Delta\to\P$ in the following way. For every
point in $\P$ we choose a homogeneous representation of the form
$(w_0:\ldots:w_n)$ with $w_0\in\{ 0,1\}$ and set
$$P\left((w_0:\ldots:w_n),(\zeta_1,\ldots,\zeta_n)\right)=
\left\{\begin{array}{ll} (1:w_1+\zeta_1:\ldots:w_n+\zeta_n),&\mbox{if}\;\;
                        w_0=1,\\
                        (0:w_1:\ldots:w_n)&\mbox{if}\;\; w_0=0.
                      \end{array}\right.$$
It is easy to see that $P$ is holomorphic in both variables and
belongs to $\Lip(1)$ with respect to each variable.
So if $G\subset\C$ is a region, $f:G\to\P$ a holomorphic curve with
Lipschitz constant $L_1$, and $\phi:G\to\Delta$ is a holomorphic map
with Lipschitz constant $L_2$ ($\Delta$ is equipped with the Euclidean metric),
then $P(f,\phi):G\to\P$ is a holomorphic curve with Lipschitz constant
$L_1+L_2$. We need two properties of $P$ for future references:
\begin{equation}
\label{inverse}
\mbox{for every}\; w\;\mbox{and every}\;\zeta\in\Delta:\quad P(P(w,\zeta),-\zeta)=w,
\end{equation}
and
\begin{equation}
\label{lips}
\dist\left( P(w,\zeta^\prime),P(w,\zeta^{\prime\prime})\right)\leq
\|\zeta^\prime-\zeta^{\prime\prime}\|.
\end{equation}

\begin{lemma}
\label{lattice} Let $E\in\C$ be a $K$-sparse set. For every $z\in\C$
we denote by $s=s(z)$ the closest point from $z$ in $E$. Then
$$\sum_{m\in E\backslash\{ s\}}\frac{1}{|z-m|^3}\leq 8K^{-3}\sum_{n=2}^\infty
(\sqrt{n}-1)^{-3}\leq 200K^{-3}.$$
and
$$\sum_{m\in E\backslash\{ s\}}\frac{1}{|z-m|^4}\leq 16K^{-4}\sum_{n=2}^\infty
(\sqrt{n}-1)^{-4}\leq 800K^{-4}.$$

\end{lemma}
{\em Proof}. If we surround every point $m\in E$ by an open disc
of radius $K/2$, centered at this point, these discs will be disjoint. We can assume that $z=0$.
Enumerate the points of $E$ is the order of increase of their distances
from the origin and
let $r_n$ be the $n$-th distance, $n=1,2\ldots$. Then $n$ points belong
to the closed disc of radius $r_n$ and the discs of radii $K/2$,
surrounding these points are disjoint and all contained in $B(0,r_n+K/2)$.
So we obtain $r_n\geq (\sqrt{n}-1)K/2$. 
\hfill$\Box$

Now we are ready to write the curve which will solve the interpolation 
problem. We apply Lemma \ref{twopoint} with $q=(1:0:\ldots:0)$ in (c)
and obtain the family of curves $g$, satisfying
(a), (b) and (c) of Lemma \ref{twopoint}. In our fixed
coordinate system $g$ has the following {\em meromorphic} homogeneous
representation:
$$g(z,a)=\left(1:g_1(z,a):\ldots:g_n(z,a)\right),$$
where $z\mapsto g_j$ are certain meromorphic functions in $\C$. 
(This representation is different from the one used in the proof of Lemma
\ref{twopoint}).
Property (c) implies that $z\mapsto g_j(z)$ are actually holomorphic
for $|z|>3$ and using (\ref{ldn}) we obtain 
\begin{equation}
\label{estim}
 \|\psi^{-1}(g(z,a))\|\leq 7|z|^{-3}<1/3,\quad |z|\geq 3
\end{equation}
and
\begin{equation}
\label{estim2}
\|\psi^{-1}(g(z,a^\prime))-\psi^{-1}(g(z,a^{\prime\prime}))\|\leq
2\| a^\prime-a^{\prime\prime}\|\,|z|^{-4},\quad |z|\geq 3,
\end{equation}
where $\psi$ was defined in (\ref{coord}).
Now we assume that a $K$-sparse set $E\subset\C$ is given with $K>25$.
An element $\a\in V^E$ is a function $E\to V,\; m\mapsto a_m$.
We define
$$f(z,\a)=
\left(1:f_1(z,\a):\ldots:f_n(z,\a)\right),$$
where
$$f_j(z,\a)=\sum_{m\in E}g_j(m-z,a_m),\quad 1\leq j\leq n.$$
Let $s\in E$ be the closest point to $z$. Then
we can single out the most important term in the previous sum:
\begin{equation}
\label{3}
f(z,\a)=P\left( g(s-z,a_s),\phi^s(z,\a)\right),
\end{equation}
where
$$\phi^s(z,\a):=\sum_{m\in E\backslash\{ s\}}\psi^{-1}(g(m-z,a_m)),$$
and $\psi^{-1}\circ g$ is well defined in view of the estimate (\ref{estim}) and
$|m-z|\geq K>25$ for $m\in E\backslash\{ s\}$.
Now we derive from (\ref{estim}), (\ref{estim2}) and Lemma \ref{lattice}
the following estimates:
\begin{equation}
\label{1}
\ds
\begin{array}{l}
\|\phi^s(z,\a))\|\leq\sum_{m\in E\backslash\{ s\}}
\|\psi^{-1}(g(m-z,a_m)\| \\
\eqnspace
\leq 7\sum_{m\in E\backslash\{ s\}}|m-z|^{-3}\leq 1400K^{-3}<1/11=\delta.
\end{array}
\end{equation}
and
\begin{equation}
\label{2}
\ds
\begin{array}{l}
\|\phi^s(z,\a^\prime)-\phi^s(z,\a^{\prime\prime})\|\leq 3\|\a^\prime-
\a^{\prime\prime}\|_\infty\sum_{m\in E\backslash\{ s\}}|m-z|^{-4}\\
\eqnspace
\leq 1600K^{-4}\|\a^\prime-\a^{\prime\prime}\|_\infty<\delta
\|\a^\prime-\a^{\prime\prime}\|_\infty,
\end{array}
\end{equation}
where $\|\a\|_\infty:=\sup_{m\in E}\| a_m\|$.
{}From (\ref{1}) and Cauchy's estimate for derivatives follows
that the spherical derivative of $\phi^s$ is uniformly bounded with
respect to $\a$. Thus by Lemma \ref{twopoint} (a) and the Lipschitz
property of $P$ we conclude that the spherical derivative of $f$ is
uniformly bounded with respect to $\a\in V^E$.

It remains to show that for every $\b\in (\P)^E$ one can find
$\a\in V^E$ such that $f(\cdot,\a)|E=\b$.
We rewrite (\ref{3}) for $z=s$ as
\begin{equation}
\label{f}
f_s(\a):=f(s,\a)=P(g(a_s),\phi_s(\a)),
\end{equation}
where $g(a):=g(0,a)$ and $\phi_s(\a):=\phi^s(s,\a)$.
{}From (\ref{1}) follows $\phi_s(\a)\in\Delta:=\{\zeta\in\C^n:\| \zeta\|<1\}$,
that is
\begin{equation}
\label{32bis}
\|\phi_s(\a)\|<\delta,\quad \a\in V^E,
\end{equation}
and from (\ref{2}) we obtain
\begin{equation}
\label{*}
\|\phi_s(\a^\prime)-\phi_s(\a^{\prime\prime})\|\leq\delta\|\a^\prime-
\a^{\prime\prime}\|_\infty,\quad \a^\prime,\a^{\prime\prime}\in V^E.
\end{equation}
%We will also use the notation $\phib:V^E\to\Delta^E$, where
%$\phib(\a),\;\a\in V^E,$
%is the function $E\to\Delta$ assuming the value $\phi_s(\a)$ at the
%point $s\in E$.

The possibility of free interpolation is now obtained from the following
version of the Inverse Function Theorem, where we put $S=\P,\;\epsilon=1,
\;L=5$ and $\delta=1/11$. 
\begin{lemma}\label{abstract}
Let $V$ and $S$ be metric spaces, $W\subset V$ a subspace with the property
$(\ref{complete})$. Let $g:V\to S$ be a map, such that $g|W$ is surjective
and in every ball $B(g(a),\delta),\; a\in W$ there exists
a right inverse to $g$, taking $g(a)$ to $a$ and with Lipschitz constant $L$.

Let $E$ be arbitrary set, and $\fb:V^E\to S^E$ be defined by 
$$f_s(\a):=f(s,\a)=P(g(a_s),\phi_s(\a)),\quad s\in E,$$
where $P: S\times\Delta\to S$ is a continuous map with
properties $(\ref{inverse})$
and $(\ref{lips})$ and $\phi_s:V^E\to\Delta$ satisfies $(\ref{32bis})$ and
$(\ref{*})$.
Assume that
\begin{equation}
\label{condition}
L\delta<\epsilon/(1+\epsilon).
\end{equation}
Then $\fb$ is surjective.
\end{lemma}

{\em Proof.} Let $\b\in S^E$. Using surjectivity of $g|W$
we find $\a^0\in W^E\subset V^E$ such that
$$g(a^0_m)=b_m,\quad m\in E.$$
Denote by $g^{-1}_m$ the right inverses to $g$, defined in
$B(b_m,\delta)$ and such that $g^{-1}_m(b_m)=a^0_m$.
We construct inductively a sequence $(\a^k),\; k=1,2\ldots$.
Assume that $\a^{k-1}$ is already defined.
Then $\|\phi_m(\a^{k-1})\|<\delta$ for $m\in E$ by (\ref{32bis}) and thus
\begin{equation}
\label{starstar}
\ds\begin{array}{l}
\dist\left( P\left( b_m,-\phi_m(\a^{k-1})\right),b_m\right)\\
\eqnspace
=\dist\left( P\left( b_m,-\phi_m(\a^{k-1})\right),P(b_m,0)\right)\leq\delta,
\end{array}
\end{equation}
where we used (\ref{inverse}) with $\zeta=0$ and (\ref{lips}).
So we can apply our right inverses $g^{-1}_m$ to
define $\a^k$ from the conditions
\begin{equation}
\label{5}
g(a^k_m)=P\left( b_m,-\phi_m(\a^{k-1})\right),\quad m\in E.
\end{equation}
Now we show by induction that
\begin{equation}
\label{33}
\|\a^k-\a^{k-1}\|_\infty\leq L^k\delta^k,\quad k=1,2\ldots.
\end{equation}
For $k=1$ this follows from (\ref{starstar}) and
the Lipschitz property of $g^{-1}_m$.
Now let $k\geq 2$. Then using (\ref{lips}) and (\ref{*}) we obtain
$$\ds\begin{array}{l}
\dist\left( P(b_m,-\phi_m(\a^{k-1})),P(b_m,-\phi_m(\a^{k-2}))\right) \\
\eqnspace
\leq \|\phi_m(\a^{k-1})-\phi_m(\a^{k-2})\|_\infty \\
\eqnspace
\leq\delta\|\a^{k-1}-\a^{k-2}\|_\infty\leq\delta L^{k-1}\delta^{k-1}=\delta^k
L^{k-1}.\end{array}$$
Applying Lipschitz property of $g^{-1}_m$ we obtain (\ref{33}).
In view of (\ref{33}) the sequence $(\a^k)$ is
a Cauchy sequence, and by (\ref{condition}) and (\ref{33})
it remains within $\epsilon$
from its original point $\a^0$, so by (\ref{complete})
it converges to some $\a\in V^E$.
In view of (\ref{5}) and continuity of $P$
$$g(a_m)=P(b_m,-\phi_m(\a)),\quad m\in E.$$
Thus by (\ref{inverse}) and (\ref{f})
$$f_m(\a)=P(g(a_m),\phi_m(\a))=P(P(b_m,-\phi_m(\a)),\phi_m(\a))=b_m,\quad
m\in E.$$
\hfill$\Box$

\vspace{.2in}

{\em eremenko@math.purdue.edu}

\end{document}